\documentclass{article}

\usepackage{amsmath,amsfonts,amssymb,latexsym,a4wide,theorem,enumerate,url}
\usepackage[utf8]{inputenc} 

{\theorembodyfont{\slshape}
\newtheorem{theorem}{Theorem}[section]

\newtheorem{lemma}[theorem]{Lemma}
\newtheorem{proposition}[theorem]{Proposition}
\newtheorem{corollary}[theorem]{Corollary}}

{\theorembodyfont{\upshape}
\newtheorem{definition}[theorem]{Definition}}

\newcommand\C{{\mathbb C}}

\newcommand\Q{{\mathbb Q}}
\newcommand\Z{{\mathbb Z}}
\newcommand\N{{\mathbb N}}
\renewcommand\P{{\mathbb P}}
\newcommand\R{{\mathbb R}}
\newcommand\HH{{\mathcal H}}
\newcommand\CC{{\mathcal C}}

\newcommand\Tau{{\mathrm T}}

\newcommand\G{{\Gamma}}
\newcommand\Lm{{\Lambda}}
\newcommand\slz{\text{SL}_2(\Z)}

\newcommand\tilX{{\widetilde{X}}}

\newcommand\eps\varepsilon
\newcommand\ph\varphi

\newcommand\qed{\hfill$\square$}
\newcommand\xalg{X^{\text{alg}}}
\newcommand\xtr{X^{\text{tr}}}
\newcommand\tr{{\text{tr}}}

\newcommand\an{\mathrm{an}}

\newcommand{\term}{\emph}

\title{Rational points on analytic  varieties} 

\author{Alexey Beshenov (Bordeaux), Margaret Bilu (Bordeaux),\\ Yuri Bilu (Bordeaux),  Purusottam Rath (Chennai)}

\date{\today}

1

\makeatletter

\renewcommand*\l@section[2]{%
  \ifnum \c@tocdepth >\z@
    \addpenalty\@secpenalty
    \addvspace{0.3em \@plus\p@}%
    \setlength\@tempdima{1.5em}%
    \begingroup
      \parindent \z@ \rightskip \@pnumwidth
      \parfillskip -\@pnumwidth
      \leavevmode \bfseries
      \advance\leftskip\@tempdima
      \hskip -\leftskip
      #1\nobreak\hfil \nobreak\hb@xt@\@pnumwidth{\hss #2}\par
    \endgroup
  \fi}
  
 \makeatother

\begin{document}

\hfuzz 4pt

\maketitle

\section*{Preface (by Yuri Bilu)}

This text is based on the expanded notes of the mini-courses on the work of  Bombieri, Pila, Wilkie, Zannier and others that  I gave in Basel in April 2011, Yaroslavl in August 2011 and Chennai in February 2012. I thank  the organizers who made it possible:
David Masser and Evelina Viada in Basel;
Fedor Bogomolov, Misha Verbitski and Alexei Zykin in Moscow;
V. Balaji, Sanoli Gun and Purusottam Rath in Chennai.

I also thank all the listeners for their precious comments. I am especially indebted to Lyosha Beshenov and Purusottam Rath, who prepared excellent notes of my courses in Yaroslavl and Chennai, respectively, and to Margaret Bilu, who attended my courses in Basel and Yaroslavl, and who wrote the present text following the notes by Lyosha and Purusottam. I only had to re-read and slightly polish the exposition at a few points, correct some inaccuracies\footnote{probably creating some new ones\dots} and add some references.

I thank Fedor Petrov and Umberto Zannier for useful discussions, Martin Klazar for detecting an error in the previous version, and Elina Wojciechowska for inspiration. 

Finally, I thank the referee for many helpful comments.

Please do not forget  that this text is based on very informal lecture notes,  and do not expect completeness and precision (in particular in the references). Please do not forget either that English is not my mother tongue. I assume full responsibility for all remaining inaccuracies, and no responsibility for the damage they may cause.

This work was supported by Bogomolov's mega-grant, ALGANT Scholarship Program, University of Basel, Erwin Schroedinger Institute (Vienna) 
and the ANR project ``HAMOT'' (ANR 2010 BLAN-0115-01).

{\footnotesize 

\tableofcontents }

\section{Introduction: Jarn\'ik and Dörge}
Let $N\in\N$ be an integer. We want to study the number of intersection points of a curve with an $N\times N$ lattice. More precisely, we want to bound the number of points of the lattice $\Lm_N = \frac{1}{N}\Z^2$ 
 lying on a curve $\G \subset [0,1)^2$ . Obviously,
$$|\G\cap \Lm_N|\leq N^2,$$
and we can find a curve $\G$ that makes this an equality, going through all vertices of the lattice in $[0,1)^2$. Suppose now that $\G$ is a graph of a function, that is, a set of the form $\{(x,y)\in [0,1)^2 \mid y = f(x)\}$. Then again, there is an easy bound
$$|\G\cap \Lm_N|\leq N,$$
and there exist graphs $\G$ realizing it. To make the problem more interesting, we are going to assume that $f$ is strictly concave or convex. The following statement is essentially due to Jarn\'ik~\cite{Ja26}. 

\begin{theorem}[Jarn\'ik, 1926]
\label{thja} 
Let $\G$ be the graph of a strictly convex or concave monotonic function ${[0,1]\to[0,1]}$. Then
\begin{equation}
\label{eja}
|\G\cap \Lm_N|\leq 3\pi^{-2/3}N^{2/3} + O(N^{1/3}\log N).
\end{equation}
\end{theorem}

The proof of Jarn\'ik's theorem will not be needed further, but it is simple and beautiful, and we sketch it below. Skipping it should not in any way harm the reader's understanding of the main content of this text. 
We thank Fedor Petrov who explained to us Jarnik's argument and his subsequent work on the topic, and Martin Klazar, who detected an error in the previous version of the proof.

It will be more convenient to prove a (formally) more general statement, by considering functions ${I\to[0,1]}$, defined on a closed sub-interval~$I$ of $[0,1]$. Thus, define 
$$
\mu(N)=\max_f|\Gamma_f\cap\Lambda_N|,
$$
where~$f$ runs over the set of strictly increasing and strictly convex functions from a closed sub-interval of $[0,1]$ to $[0,1]$, and~$\Gamma_f$ denotes the graph of~$f$. Further, let~$I$ be a shortest sub-interval of $[0,1]$ such that there exists ${f:I\to[0,1]}$ as above with ${|\Gamma_f\cap\Lambda_N|=\mu(N)}$. We fix these~$I$ and~$f$ and write ${\Gamma=\Gamma_f}$.

Call a vector ${(x,y)\in \Z^2}$ \term{primitive} if ${d=\gcd(x,y)=1}$; in general, we call $(x/d,y/d)$ the \term{primitivization} of $(x,y)$. Let~$S_X$ denotes the set of primitive vectors with positive coordinates whose sum is  bounded by~$X$:
$$
S_X=\{(x,y)\in \Z^2:\ x,y>0,\ \gcd(x,y)=1,\ x+y\le X\}
$$ 
Then
\begin{equation}
\label{esx}
|S_X|=\frac3{\pi^2}X^2+O(X\log X), \qquad \sum_{(x,y)\in S_X}x=\frac1{\pi^2}X^3+O(X^2\log X).
\end{equation}
(These relations are pretty standard. See, for instance, \cite[Section~18.5]{HW08} for the proof of a similar statement; both relations~\eqref{esx} can be proved analogously.)

Let ${P_0, \ldots, P_n}$ be the lattice points on~$\Gamma$ (so that ${n+1=\mu(N)}$) ordered by the increasing first coordinates. 
The vectors ${(x_i,y_i)=N(P_i-P_{i-1})}$ are pairwise distinct (because the function is convex) and have positive coordinates (because it is increasing). We may also assume them to be primitive; if one of them is not, then, replacing it by the primitivization, we may find a new (strictly increasing and strictly convex) function  with graph having  ${n+1}$ lattice points, but defined on a shorter interval, 
contradicting our choice of~$I$.

Let $X$ be the largest real number such that ${n+1\ge |S_X|}$. Then
$$
\begin{gathered}
|\Gamma\cap\Lambda_N|=n+1=\frac3{\pi^2}X^2+O(X\log X), \\
2N\ge\sum_i(x_i+y_i)\ge 2\sum_{(x,y)\in S_X}x=\frac2{\pi^2}X^3+O(X^2\log X),
\end{gathered}
$$
whence the result.

\bigskip

It follows from the proof that bound~\eqref{eja} is optimal in the sense that for each~$N$ we can find a curve $\G_N$ (satisfying the assumptions of Theorem~\ref{thja})
such that 
$$
|\G_N\cap\Lm_N| \ge 3\pi^{-2/3}N^{2/3}(1+o(1)).
$$
Here the curve~$\Gamma_N$ depends on~$N$. One may wonder whether there  exists a ``universal Jarn\'ik curve'' $\G$ satisfying ${|\G_N\cap\Lm_N| >\kappa N^{2/3}}$ for infinitely many  $N$, where~$\kappa$ is a positive real number, not necessarily equal to $3\pi^{-2/3}$. This question, which  is attributed to J.-M. Deshouillers and A.~Plagne~\cite{Pl99}, was  recently answered negatively by  F.~Petrov~\cite{Pe06}, who showed that for any fixed $\G$ we have
$$|\G\cap\Lm_N| = o(N^{2/3})\quad \text{as}\quad N\longrightarrow \infty.$$

Now we  give a totally different  proof of a weaker version of Jarn\'ik's theorem. This proof is based on an idea of Dörge~\cite{Doe27} (1927) (see also \cite[Section 9.2]{La83}) that will turn out to be very fruitful.

 \begin{theorem}(Weak Jarn\'ik) Let $I$ be a compact interval and $\G$ be the graph $\{(x,y) \mid y = f(x)\}$ of a strictly convex function $f\in\CC^2(I)$. Then
 $$|\G\cap\Lm_N|\leq c(I,f)\,N^{\frac23}.$$
 \end{theorem}
Note that this is weaker than Jarn\'ik's theorem because $f$ is assumed to be $\CC^2$ and because the constant in the bound depends on $f$. 

The proof uses the following generalization of Lagrange's mean value theorem:

\begin{proposition}(Schwarz mean value theorem) Let~$I$ be an interval, ${f\in\CC^n(I)}$, and ${x_0, \ldots, x_n \in I}$. Then there is $\tau\in I$ such that
$$ \left|\begin{array}{ccccc}
1 & x_0 & \cdots & x_0^{n-1} & f(x_0) \\
\vdots & \vdots & \ddots & \vdots & \vdots \\
1 & x_n & \cdots & x_n^{n-1} & f(x_n) \\
\end{array}\right| = \frac{f^{(n)} (\tau)}{n!} \, V(x_0, \ldots, x_n).$$
\end{proposition}
Here $V(x_0,\ldots,x_n)$ denotes the Vandermonde determinant
$$ \left|\begin{array}{ccccc}
1 & x_0 & \cdots & x_0^{n-1} & x_0^n \\
\vdots & \vdots & \ddots & \vdots & \vdots \\
1 & x_n & \cdots & x_n^{n-1} & x_n^n \\
\end{array}\right|. $$

\noindent\textit{Proof of Schwarz mean value theorem.} Let $g(t) = a_nt^n + \cdots + a_1t + a_0$ be the Lagrange interpolation polynomial for $f$ at the points $x_0,\ldots, x_n$. Then $(f-g)(x_i) = 0$ for all $i = 0,\ldots n$, and therefore by a ``generalized Rolle theorem'' there is $\tau\in I$ such that $f^{(n)}(\tau) = g^{(n)}(\tau)$.  Moreover, $g^{(n)}(\tau) = n!\,a_n$. Noting that 
$$V(x_0,\ldots,x_n)\,a_n  =  \left|\begin{array}{ccccc}
1 & x_0 & \cdots & x_0^{n-1} & f(x_0) \\
\vdots & \vdots & \ddots & \vdots & \vdots \\
1 & x_n & \cdots & x_n^{n-1} & f(x_n) \\
\end{array}\right|$$
then yields the result.\qed

\bigskip

\noindent\textit{Proof of the weak Jarn\'ik theorem.} Let us pick three points ${P_i = (x_i,y_i)}$, ${i = 0,1,2}$ on $\G$. Because of the strict convexity of $f$, they cannot lie on one straight line, so
$$ \Delta = \left|\begin{array}{ccc}
1 & x_0 & y_0 \\
1 & x_1 & y_1 \\
1 & x_2 & y_2 \\
\end{array}\right| \ne 0. $$
By the Schwarz mean value theorem, there exists ${\tau\in I}$ such that
$$ \Delta = \left|\begin{array}{ccc}
1 & x_0 & f(x_0) \\
1 & x_1 & f(x_1) \\
1 & x_2 & f(x_2) \\
\end{array}\right| = \frac{f^{\prime\prime}(\tau)}{2!}\,V(x_0,x_1,x_2). $$
Since $I$ is compact, this gives us an upper bound
$$|\Delta|\leq c(I,f)\max_{1\leq i < j \leq 3} |x_i-x_j|^3.$$
Assume now that ${P_0,P_1,P_2\in \G \cap\Lm_N}$, that is, the $P_i$ are rational points whose coordinates have denominators dividing $N$. Then $\Delta$ is a non-zero rational number with denominator at most $N^2$, so $$|\Delta|\geq \frac{1}{N^2}.$$ This, combined with the above upper bound, gives us
$$\max_{1\leq i < j \leq 3} |x_i-x_j| \geq c\,N^{-2/3}.$$ We have thus proved that an interval containing the projections of three points lying in $\G \cap\Lm_N$ has length bounded from below by $c N^{-2/3}$. Cutting $I$ into small intervals of length $\frac{c}{2}N^{-{2}/{3}}$, we have at most two points projecting to each such small interval. Denoting by $\ell(I)$ the length of $I$, there are at most $ \ell (I)\frac{2}{c}N^{{2}/{3}}$ such intervals, so the result follows.\qed

In fact, we have proved the following
\begin{proposition}Let $I$ be a compact interval and $\G$ the graph $\{(x,y) \mid y = f(x)\}$ of a  function $f\in\CC^2(I)$. There is a constant $c$ depending only on~$I$ and~$f$, having the following property. If for the points
$$P_0,\ldots,P_k \in\G\cap\Lm_N,\ \ \ P_i = (x_i,y_i)$$
the coordinates $x_1,\ldots,x_k$ lie in an interval of length smaller than $c\,N^{-{2}/{3}}$, then $P_0,\ldots,P_k$ lie on a straight line.\end{proposition}

In particular, the set $\G\cap\Lm_N$ is covered by $O_f(N^{2/3})$ straight lines. Bombieri and Pila~\cite{BP89} came up with a crucial idea:  consider not only straight lines, but also more general algebraic curves. We discuss it in the next section. 

\section{The theorem of Bombieri and Pila}
Generalizing Dörge's ideas, Bombieri and Pila~\cite{BP89} proved  that a plane analytic compact curve cannot contain many rational points with a given denominator.

In the sequel by an \emph{irreducible plane curve of degree~$d$} we mean a subset of $\R^2$ consisting of points $(x,y)$ satisfying $F(x,y)=0$, where $F(X,Y)$ is an irreducible real polynomial of degree~$d$. A \emph{segment} of such a curve is its connected compact subset.

\begin{theorem}\label{bombieri}(Bombieri, Pila, 1989) Let $\G\subset\R^2$ be a real analytic plane compact curve. 
\begin{enumerate}\item Assume that $\G$ is transcendental.  Then for any $\eps>0$ there is a constant $C(\G,\eps)$ such that
$$|\G\cap\Lm_N|\leq C(\G,\eps)\,N^{\eps},$$
where  $\Lm_N = \frac{1}{N}\Z^2$.
\item Assume that $\G$ is a segment of an irreducible plane algebraic curve of degree $d$. Then for any $\eps>0$ there is a constant $C(\G,\eps)$ such that
$$|\G\cap\Lm_N|\leq C(\G,\eps)\,N^{\frac{1}{d}+\eps}.$$
\end{enumerate}
\end{theorem}
\paragraph{Remarks}\begin{enumerate} \item The example of a curve given by $y =x^d$ shows that the exponent in the second part of the theorem cannot be made smaller than ${1}/{d}$.
\item In the second part of the theorem one can replace $C(\G,\eps)$ by $C(d,\eps)$, if one assumes in addition that $\G\subset[0,1]^2$.
\end{enumerate}

We are going to use the following generalization of the mean value theorem:

\begin{proposition} (Generalized Schwarz mean value theorem) Let~$I$ be an interval, ${f_0,\ldots,f_n\in \CC^n(I)}$, and ${x_0,\ldots,x_n\in I}$. Then for all $i,j$ such that ${0\leq i,j \leq n}$ there is $\tau_{ij}\in I$ such that
$$ \det [f_i (x_j)]_{i,j} = \det \left[\frac{f_i^{(j)} (\tau_{ij})}{j!}\right]_{ij} \cdot V (x_0, \ldots, x_n). $$
\end{proposition}
 Here is the proof  for $n=1$: by the Lagrange mean value theorem, there are ${\tau_{01},\tau_{11}\in I}$ such that ${f_i(x_1) - f_i(x_0) = f'_i(\tau_{i1})\,(x_1 - x_0)}$ for $i=0,1$. Thus,
$$ \left|\begin{array}{cc}
f_0 (x_0) & f_0 (x_1) \\
f_1 (x_0) & f_1 (x_1) \\
\end{array}\right| =
\left|\begin{array}{cc}
f_0 (x_0) & f_0 (x_1) - f_0 (x_0) \\
f_1 (x_0) & f_1 (x_1) - f_1 (x_0) \\
\end{array}\right| =
\left|\begin{array}{cc}
f_0 (x_0) & f_0^\prime (\tau_{01}) \\
f_1 (x_0) & f_1^\prime (\tau_{11}) \\
\end{array}\right|\,(x_1-x_0), $$
and we can choose $\tau_{00} = \tau_{10} = x_0$.

The proof of the general case is very similar, but with Lagrange replaced by Schwarz, see \cite{vieux}. 
\begin{corollary}\label{cor} If the interval $I$ is compact, there is a constant $c$ depending only on $f_0,\ldots,f_n$ such that
$$\left|\det [f_i (x_j)]_{i,j}\right|\leq c\,|V(x_0,\ldots,x_n)|.$$
\end{corollary}
\begin{definition} Let $d,e\in\N$. We call a \term{$(d,e)$-curve} a (possibly reducible) plane algebraic curve with equation of the form $P(x,y) = 0$ where $P\in\R[x,y]$ is a non-zero polynomial such that $\deg_x P < d$ and $\deg_y P < e$.
\end{definition}
For example, a $(2,2)$-curve is of the form $a + bx + cy + hxy = 0$.

\paragraph{Facts about $(d,e)$-curves}\begin{enumerate}\item Any collection of $de-1$ points lies on a $(d,e)$-curve.
\item A collection of $de$ points $P_k = (x_k,y_k)$, $k = 1,\ldots,de$, lies on a $(d,e)$-curve if and only if 
$$ \Delta := \det [x_k^i\,y_k^j]_{\substack{(i,j) = (0,0), \ldots, (d-1,e-1)\\ k = 1, \ldots, d\,e}} = 0. $$
Note that the rows in the above determinant are numbered by the bi-index $(i,j)$. For example, for $(2,2)$-curves this condition for four points $P_1,\ldots, P_4$ takes the form
$$ \left|\begin{array}{cccc}
1 & 1 & 1 & 1 \\
x_1 & x_2 & x_3 & x_4 \\
y_1 &  y_2& y_3 & y_4 \\
x_1\,y_1 & x_2\,y_2 & x_3\,y_3 & x_4\,y_4 \\
\end{array}\right| = 0. $$
\item If $P_1,\ldots,P_s$ do not lie on a $(d,e)$-curve (in particular $s\geq de$  by the above), then one can select $de$ points among them not lying on a $(d,e)$-curve.
\end{enumerate}
Let $\G$ be the graph $y = f(x)$ of a function $f\in \CC^{de-1}(I)$ where $I$ is an interval, and consider $P_1,\ldots,P_{de}\in\G\cap\Lm_N$ not lying on a $(d,e)$-curve. Then, by the above, writing $g_{i,j}(x) = x^if(x)^j$ we have
$$ \Delta = \det [g_{ij}(x_k)]_{\substack{(i,j) = (0,0), \ldots, (d-1,e-1)\\ k = 1, \ldots, d\,e}} \ne 0, $$
where $P_i = (x_i,f(x_i))$. Using Corollary \ref{cor} for the functions $g_{ij}$, we get
$$|\Delta|\leq c(f)\,|V(x_1,\ldots,x_{d\,e})|\leq c\left(\max_{1\leq k<l\leq d\,e}|x_k-x_l|\right)^{\frac{de(de-1)}{2}},$$
where the constant $c$ depends only on~$I$ and~$f$. On the other hand, the coordinates of $P_i$'s are rational numbers with denominator dividing $N$. By assumption, the determinant $\Delta$ evaluates to a non-zero rational number, and the entries in its $(i,j)$-th column have a common denominator dividing $N^{i+j}$. Thus the denominator of $\Delta$ is bounded by 
$$\prod_{(i,j)} N^{i+j} = N^{\sum_{(i,j)} i+j} = N^{\frac{de(d+e-2)}{2}}$$
and we have
$$|\Delta|\geq N^{-\frac{de(d+e-2)}{2}}.$$
We can therefore conclude that if points $P_1, \ldots, P_{de}\in\G\cap\Lm_N$ do not lie on a $(d,e)$-curve, then
$$\max |x_k - x_l| \ge c(f) \cdot N^{- \frac{e+d-2}{ed-1}} = c(f) \cdot N^{-\delta},$$
where $\delta = \frac{d+e-2}{de-1}.$
This means that if $de$ points of $\G\cap\Lm_N$ project to an interval of length $\frac{c}{2}N^{\delta}$, they necessarily lie on a $(d,e)$-curve. Since we can cut our interval $I$ into $c(f,d,e,I)\,N^\delta$ such intervals, we have the following.

{\sloppy

\begin{lemma}(Main lemma of Bombieri--Pila) Let~$\G$ be the graph ${y = f(x)}$ of a function ${f\in \CC^{de-1}(I)}$. Then the set $\G\cap\Lm_N$ is covered by $cN^\delta$ \quad  $(d,e)$-curves, where~$c$ depends on~$I$,~$f$,~$d$,~$e$ and  ${\delta = \frac{d+e-2}{de-1}}$. 
\end{lemma}
Note that the approach of Bombieri and Pila is an exact analogue of the proof of the weak Jarn\'ik theorem we gave above, with straight lines replaced by more general $(d,e)$-curves. 

}

\bigskip
\noindent\textit{Proof of Theorem \ref{bombieri}.}\ We prove the two statements of the theorem separately:
 \begin{enumerate}\item Take $\G$ transcendental, analytic and compact as in the statement of the theorem, and $\eps>0$. Choose moreover $d,e$ such that ${\frac{d+e-2}{de-1}<\eps}$. We can break up $\G$ into a finite number of pieces of the form $y = f(x)$ or $x = f(y)$, to each of which we apply the main lemma. Thus $\G \cap\Lm_N$ is covered by $cN^\eps$ \quad $(d,e)$-curves, where ${c= c(\G,\eps)}$. On the other hand, since $\G$ is transcendental, the size of the intersection of $\G$ with a $(d,e)$-curve is bounded uniformly by a constant depending only on $\G$, $d$, $e$.

\item Assume now that $\G$ is a segment of an irreducible algebraic curve of degree $d$. By applying a linear change of variables, we may assume that it is defined by the polynomial equation ${F(x,y)=0}$ where ${\deg_xF=d}$. 
Choose $e$ big enough such that
$\frac{d+e-2}{de-1}<\frac{1}{d} + \eps.$ Then, applying the main lemma, $\G\cap\Lm_N$ is covered by $cN^{1/d+\eps}$ \quad  $(d,e)$-curves, where ${c=c(\G,\eps)}$. On the other hand, since ${\deg_xF=d}$, our curve~$\Gamma$ cannot be contained in a $(d,e)$-curve. By the Bézout theorem, the size of the intersection of $\G$ with a $(d,e)$-curve is bounded by $c(d,e)$, and we are done. \qed
\end{enumerate}

\section{Counting rational points with bounded height}
In this section, instead of counting lattice points on $\G$, we are going to count points of bounded height.

\begin{definition} Let $\alpha = \frac{p}{q}\in\Q$ with $p$ and $q$ coprime. Define the height of $\alpha$ by
$$H(\alpha) = \max\{|p|,|q|\}.$$
\end{definition}
Denoting $\G(\Q,N) = \{(x,y)\in\G\cap\Q^2,\ H(x),H(y)\leq N\}$, we can state the following stronger version of Theorem \ref{bombieri}:

\begin{theorem}
(Bombieri--Pila theorem for rational points) 
\label{tbprat}
Let $\G$ be a plane compact analytic curve.\begin{enumerate} \item If $\G$ is transcendental, then for any $\eps>0$ there is a constant $c(\G,\eps)$ such that
$$|\G(\Q,N)|\leq c(\G,\eps)\,N^{\eps}.$$
\item If $\G$ is algebraic of degree $d$, then for any $\eps>0$ there is a constant $c(\G,\eps)$ such that
$$|\G(\Q,N)|\leq c(\G,\eps)\,N^{\frac{2}{d}+\eps}.$$
\end{enumerate}
\end{theorem}
Again, the second part cannot be improved, as shown by the example $y = x^d$.

\bigskip

\noindent \textit{Proof.} To get this result, we just need to modify a little the proof of Theorem \ref{bombieri}, more precisely, the proof of the Main Lemma.  Consider points
$$P_1,\ldots,P_{de}\in \G(\Q,N),\ \ \ P_k = (x_k,y_k)$$
not lying on a $(d,e)$-curve, where $\G$ is a graph $y=f(x)$. We still have the bound
$$|\Delta| \le c(f)\left(\max_{1\leq k<l\leq d\,e}|x_k-x_l|\right)^{\frac{de(de-1)}{2}}$$
obtained directly from the general mean value theorem where
$$ \Delta = \det [x_k^iy_k^j]_{\substack{(i,j) = (0,0), \ldots, (d-1,e-1)\\ k = 1, \ldots, d\,e}}\ne 0. $$
In fact, the only thing that needs to be changed is the computation of the lower bound for $\Delta$.
The difference is that since now our condition is on the height of $x_k$'s and $y_k$'s, we have less information on their denominators, and in particular we can no longer use $N$ as a common denominator. Write $(x_k,y_k) = \left(\frac{\cdot}{M_k},\frac{\cdot}{N_k}\right)$, with $0<M_k,N_k \le N$. Then the common denominator in column $k$ is at most $M_k^{d-1}N_k^{e-1}\le N^{d+e-2}$, and $\Delta$ has denominator at most $N^{de(d+e-2)}$, which gives us
$$|\Delta|\ge N^{-de(d+e-2)}. $$
The rest of the proof works in exactly the same way. Note that in the proof of the first version of the Main Lemma, the exponent was $\frac{de(d+e-2)}{2}$, that is, we lost the factor $\frac{1}{2}$. Therefore, the new version of the Main Lemma is the following:
\begin{lemma}(Second Main Lemma) Let $I$ be an interval and $\G$ a graph $y = f(x)$ of a function $f\in\CC^{de-1}(I)$. Then $\G(\Q,N)$ is covered by $c\,N^{\delta}$ \quad $(d,e)$-curves, where $c$ depends only on $I$, $f$, $d$ and $e$, and $\delta = \frac{2(d+e-2)}{de-1}$.
\end{lemma}

\section{Generalization to higher dimensions}

\label{shigher}

Let $X\subset\R^m$ be a compact real-analytic transcendental manifold. We denote by $X(\Q,N)$ the set of rational points of $X$ whose coordinates are all of height smaller than $N$. To generalize the previous results, we would like to get something like
$$|X(\Q,N)|\le c(X,\eps) \, N^{\eps}.$$
However this is easily seen to be false, as the example of the surface $z = x^y$ where $1\leq x,y\leq 2$ shows. Indeed, the problem that arises in the higher dimensional case is the fact that a transcendental manifold can very well contain whole pieces of algebraic curves, which can have many rational points. Let us therefore define $\xalg$ to be the union of all segments of algebraic curves inside $X$, and $\xtr = X\setminus \xalg$. The above statement then turns out to be true when replacing $X$ by $\xtr$:

\begin{theorem}\label{higher}(Pila, Wilkie) For every $\eps>0$ there is a constant $c(X,\eps)$ such that
$$\xtr(\Q,N)\le c(X,\eps)N^{\eps}.$$
\end{theorem}
We  sketch now the proof in the case where $X$ is a (two-dimensional) surface inside $\R^3$. In this special case, Theorem \ref{higher} was proved by Pila in \cite{Pi05}. 

 \begin{theorem}\label{pila05}(Pila, 2005) Let $X\subset\R^3$ be an analytic surface.  For every $\eps>0$, there is a constant $c(X,\eps)$ such that
$$\xtr(\Q,N)\le c(X,\eps)\,N^{\eps}.$$
 \end{theorem} 
 
In fact, Pila proved this result for subanalytic surfaces, projections of analytic surfaces. See \cite{BM88} for a tutorial on these topics. 

We start from an analog of the Main Lemma. Instead of algebraic curves of bounded degree, as in the one-dimensional case, we now use algebraic surfaces of bounded degree: for a positive integer $d$, a \textsl{$d$-surface}\footnote{To make it compatible with the $(d,e)$-curves, we should have probably said \textsl{$(d,d,d)$-surface}, but this is too lengthy.} is defined by $P(x,y,z) = 0$ where
$\deg_x P,\, \deg_y P,\, \deg_z P < d$. The Main Lemma asserts that rational points on a sufficiently smooth compact surface lie on a few $d$-surfaces. 

\begin{proposition}(Main Lemma in dimension~$2$)
Let $\eps >0$. There are integers $d$ and $D$ depending on $\eps$, such that if $X$ is a compact $\CC^D$-surface in $\R^3$, then $X(\Q,N)$ is covered by $c\,N^\eps$ \quad $d$-surfaces, where ${c=c(X,\eps)}$.
\end{proposition}

The proof goes along the same lines as in dimension~$1$. Fix~$d$ and~$D$, to be specified later. Points $P_1,\ldots,P_s\in \R^3$ belong to a $d$-surface if a certain determinant $\Delta$ vanishes. We may assume that~$X$ is given by ${z = f(x,y)}$,  where $f$ is a function   defined on some compact domain in $\R^2$ having bounded continuous derivatives of order up to~$D$. If  $P_1,\ldots,P_s\in X$ project to a small square on the $x,y$-plane, and~$D$ is large enough (in terms of~$d$), then $\Delta$ can be bounded from above using some analogues of the Mean Value Theorem. On the other hand, if $P_1,\ldots,P_s$ are rational points of height at most~$N$ and $\Delta\ne 0$, then $\Delta$ is bounded from below by some negative power of~$N$. When the square is small enough, the upper bound contradicts the lower bound, which means that $\Delta$ must be $0$. In other words, if the points project to a sufficiently small square, then they must belong to a $d$-surface. Selecting~$d$ suitably large, we see that ``sufficiently small square'' means ``square with   side length $c(X,\eps)N^{-\eps}$\;''. Now subdividing the domain of definition of the function~$f$ into small squares, we complete the proof. For the details see \cite[Appendix A, Lemma A.3]{Za12}. 

\bigskip

Having this, let us see what kind of intersection $X$ can have with a $d$-surface. Since $X$ is transcendental, such an intersection is of dimension~$1$; by compactness, it must be a finite union of irreducible analytic curves. The number of irreducible components in the intersection is controlled by the following result due to Gabrielov~\cite{Ga68}. 

\begin{theorem}(Gabrielov) Let $V,W$ be complex semi-analytic manifolds with $W$ compact, and $\pi:V\to W$ an analytic map. There exists a constant $c(\pi)$ depending only on $\pi$ such that for all $\omega \in W$ the number of connected components of $\pi^{-1}(\omega)$ is bounded by $c(\pi)$.
\end{theorem}

Applying this theorem (more precisely, an analogous statement about irreducible components) with~$W$ being the space of all non-zero polynomials defining $d$-surfaces (which is $\P^{d^3-1}(\R)$, hence compact), we bound the number of components in the intersection of~$X$ with a $d$-surface. Now let us consider one such component. It is a curve, algebraic or transcendental. If it is algebraic, it is contained in $\xalg$ and we don't have to worry about it. If it is transcendental, then we reduced the problem to the one-dimensional case. Projecting it to one of the coordinate planes, we obtain a plane transcendental curve, and to the latter we may apply Theorem~\ref{tbprat}. 

The argument above can be briefly summarized as follows: rational points of bounded height from~$X^\tr$ belong to ``few'' transcendental curves, and each of the latter has ``few'' rational points of bounded height by Theorem~\ref{tbprat}. 

However, there is an important difficulty: Theorem~\ref{tbprat}  involves a constant depending on the curve. Thus we need to have a sort of ``uniform'' version of this theorem, which would follow if we get a similar uniform version of the one-dimensional Main Lemma.

Let us therefore for a moment go back to the one-dimensional case. Let~$\Gamma$ be the graph  ${y=f(x)}$, where~$f$ is defined on some interval~$I$ and has sufficiently many bounded derivatives therein. The Main Lemma asserts that  
the set  $\G(\Q,N)$ is contained inside a union of $cN^{\eps}$ \quad $(d,d)$-curves, where $c$ depends on~$I$,~$f$,~$d$ and~$\eps$.

Now assume that~$f$ runs through a continuous family of functions ${\{f_\tau:I_\tau\to \R: \tau\in \Tau\}}$ with a compact base~$\Tau$. We will manage to complete the argument above if we  bound the constant~$c$ uniformly in~$\tau$. This would have been possible if the constant depended only on the length $|I|$ of the interval and the sup-norm $\|f\|$: both are bounded in a continuous family with compact base. Unfortunately, the constant depends also on the  sup-norms of the derivatives $\|f^{(k)}\|$ with ${k = 0,\ldots,d^2-1}$, and the derivatives are not bounded in a  continuous family with compact base (for example, take ${\Tau=[0,1]}$, ${I_\tau= [0,1]}$ and ${f_\tau(x)=\sqrt{x+\tau}}$). 

One can take care of the first derivative $f'$ by noticing that on the part of the interval where $|f'|\ge1$ our curve can be written as ${x=g(y)}$ with ${|g'|\le 1}$. However, one cannot deal like this with higher order derivatives.

Instead, Pila uses the fact that derivatives can be large only on short intervals. The following lemma is easy to prove.

\begin{lemma} Assume that $f$ is of class $\mathcal{C}^2$ on  $I=[a,b]$, and $f,f',f''$ do not vanish on $I$.  Put $J = [a+\delta,b-\delta]$ where $0 < \delta < \frac{b-a}{2}$. Then
${\|f'\vert_J\|< \delta^{-1}\|f\|}$.
\end{lemma}

Indeed, assume, for instance, that $f,f',f''>0$ on $I$. Then for ${x\in J}$  we have ${f(x)> 0}$ and ${f'(u)> f'(x)>0}$ for ${u\in [x,b]}$. Hence
$$
f(b) = f(x)+\int_x^b f'(u)\,du > f'(x)\,(b-x)\ge f'(x)\,\delta,
$$
whence ${0<f'(x)< \delta^{-1}f(b)= \delta^{-1}\|f\|}$, proving the lemma in this case. The other cases are treated similarly. 

By induction, one shows that if ${f\in \mathcal{C}^{n+1}(I)}$ and $f,f',\ldots,f^{(n+1)}$ do not vanish on $I$, then ${\|f^{(n)}\vert_J\|\le (n/\delta)^n\,\|f\|}$. 

Thus high order derivatives may grow uncontrollably only near the points where one of them vanishes.

Now return to our curve $y=f(x)$. Put ${n=d^2-1}$ and assume that $f,f',\ldots, f^{n+1}$ have $k$ roots altogether. Fix a small positive number~$\delta$ and throw away from our interval $I$ all the $\delta$-neighborhoods  of these roots. On the remaining part of~$I$ the sup-norms of derivatives of order up to~$n$ are controlled by the sup-norm of~$f$. As for the part thrown away, it consists of~$k$ tiny intervals, and Pila applies to them an ingenious re-scaling argument, going back to the original article with Bombieri. The details are quite intricate and cannot appear here: the reader may consult Pila's article \cite{Pi05} or the exposition in \cite[Appendix A]{Za12}.

This reasoning implies a new version of the Main Lemma, with constant depending only on $|I|, \|f\|$ and the number of zeros of derivatives of order up to $n$. For the curves occurring in the proof of Theorem~\ref{pila05} the first two of these parameters are estimated immediately, just from the compactness of the space of $d$-surfaces. As for the number of zeros of derivatives, this can be estimated using Gabrielov's theorem. 

This is, in general terms, how Theorem~\ref{pila05} was proved. 

\section{The Manin--Mumford conjecture}

In this section we show how Pila's techniques apply to the famous problem of Manin--Mumford. We start with some generalities about complex abelian varieties. 

Let $A$ be an abelian variety over $\C$ of dimension $g$. Recall that there is a complex analytic group isomorphism 
$$A(\C)\cong \C^g/\Lm,$$
where $A(\C)$ is the group of complex points of $A$ and $\Lm$ is a lattice inside $\C^g$. It follows  that
$$A[N]\cong \left(\Z/N\Z\right)^{2g},$$
where $A[N] = \{x\in A \mid Nx = 0\}$  are the $N$-torsion points of $A$. We will denote the group of all torsion points of $A$ by $A^{\text{tors}} = \cup_{N\ge 1} A[N].$
 In its simplest form, the Manin--Mumford conjecture can be stated as follows:
\begin{theorem}
\label{tmmsi}
(Manin--Mumford conjecture, simplest form) Let $X$ be an algebraic curve on an abelian variety~$A$. If $X$ is not an elliptic curve, then $|X\cap A^{\text{tors}}|$ is finite.
\end{theorem}
If $X\subset A$ is an elliptic curve passing through a torsion point of~$A$, then the above set is equal to the set of torsion points of $X$, so is infinite: this explains why we need to exclude elliptic curves. 

Theorem~\ref{tmmsi} was proved by Raynaud~\cite{Ra83} in 1983. To state a more general form of the Manin--Mumford conjecture, we need the following definition.

\begin{definition} A subvariety $B\subset A$ is called a \term{torsion subvariety} if it is a translate of an abelian subvariety by a torsion point.  That is,  there exist  a torsion point $b\in A$ and an abelian subvariety $B_0\subset A$ such that $B = b + B_0$.
\end{definition}
\begin{theorem} (Manin--Mumford conjecture, general form) Let $X \subset A$ be a complex algebraic variety. Then $X$ has only finitely many maximal torsion subvarieties.
\end{theorem}
This was proved by Raynaud~\cite{Ra83a} in 1983 as well. After the work of Raynaud a number of other proofs emerged. 

Recently a new proof was suggested by Pila and Zannier (2006). They do prove the general conjecture, but we will discuss their argument only  in the simplest case where $X$ is an algebraic curve.  Let us consider the inverse image of $X$ under the the complex analytic uniformization map for $A$
\begin{eqnarray*}
\pi:\mathbb{C}^g & \longrightarrow& A,\\
\frac{1}{N}\Lm & \longrightarrow &A[N],\\
\tilX = \pi^{-1} (X) & \longrightarrow & X.
\end{eqnarray*}
Let ${\Delta\subset \C^g}$ be a fundamental domain of the lattice $\Lm$. Then $\pi\vert_{\Delta}:\Delta\longrightarrow A$ is surjective, and any point inside $A[N]\cap X$ comes from a point in $\tilX\cap\frac{1}{N}\Lm\cap\Delta$, so we have 
$$\bigl|X\cap A[N]\bigr|\le\left|\tilX\cap\frac{1}{N}\Lm\cap\Delta\right|.$$
Thus we have reduced the problem to rational points on the complex analytic curve (and real analytic surface) $\tilX$. Our goal is to apply Theorem \ref{higher}, for which we have to determine~$\tilX^{\mathrm{alg}}$. This will be done in two steps:

\paragraph{Claim 1:} Each irreducible component of~$\tilX$ is transcendental.

\bigskip

First of all, let us observe  the following property of plane algebraic curves.

\begin{lemma}
Let ${C\subset \R^2}$ be an irreducible plane algebraic curve and let $\Lambda$ be a lattice in $\mathbb{R}^2$. Then either $C$ is a straight line and the intersection ${C\cap\Lambda}$ is infinite, or ${C+\Lambda}$ is dense in $\R^2$ (in the real topology).
\end{lemma}

We do not prove this lemma, but it is not very difficult.  Since $\tilX$ is invariant under translations by elements of the lattice, but not dense, the lemma implies that its only possible algebraic components are straight lines having infinite intersection with the lattice. But having such a component would mean that $X$ has an elliptic curve as a component, a contradiction. This proves Claim 1.

\paragraph{Claim 2:} $\tilX$ cannot contain a segment of a real algebraic curve.

\bigskip

Indeed, by analytic continuation this would give us a complex algebraic curve inside $\tilX$, contradicting Claim 1.

\bigskip 

Thus $\tilX^{\mathrm{alg}} = \varnothing$. We are therefore ready to apply the Pila--Wilkie theorem to get an upper bound
$$\bigl|X\cap A[N]\bigr|\le\left|\tilX\cap\frac{1}{N}\Lm\cap\Delta\right|\le c(X,\eps)\,N^{\eps}.$$
(Notice that while $\tilX$ can a priori have infinitely many components, only finitely many of them may intersect $\Delta$.)

For the lower bound, we may assume without loss of generality that $A$ and $X$ are defined over a number field $K$. Then if $P\in X\cap A[N]$, all conjugates of $P$ over~$K$ lie in $X\cap A[N]$ as well, which implies $|X\cap A[N]|\ge [K(P):K]$ where $K(P)$ is the extension of $K$ generated by the coordinates of $P$. The lower bound comes from the fact that torsion points generate extensions of large degree, an old result of Masser~\cite{Ma84}. 
\begin{theorem} (Masser) Let $A$ be an abelian variety of dimension $g$ over a number field $K$. Then there exist $\delta = \delta(g)>0$ and $c = c(A,K)>0$ such that if $P\in A[N]$ is of exact order $N$ then
$$[K(P):K]\ge cN^{\delta}.$$
\end{theorem}
Taking $\eps = \frac{\delta}{2}$ shows that for sufficiently big $N$, the set $X\cap A[N]$ must be empty.

\section{Definable sets in o-minimal structures}

In fact, in \cite{PW06} Pila  and Wilkie obtained a much more general result than Theorem~\ref{higher}: they estimated the number of rational points of bounded degree on    arbitrary  definable sets in o-minimal structures over $\R$.

\subsection{O-minimal structures}

For our purposes the following definitions will be sufficient; the interested reader might want to consult the main reference on the subject, \cite{dries}.

\begin{definition}We define a \term{structure} on a set $R$ to be a sequence
$$\mathcal{S}_0, \mathcal{S}_1, \mathcal{S}_2, \mathcal{S}_3, \ldots,$$
where $\mathcal{S}_m\subset \mathcal{P}(R^m)$ is a collection of subsets of $R^m$.
A set $A \subseteq R^m$ such that $A\in \mathcal{S}_m$ is called \term{definable} in the structure.
Further we ask that

\begin{enumerate}[(S1)]\item Each $\mathcal{S}_m$ is closed under finite union and complement (and therefore  under finite intersection as well): if $A, B \in \mathcal{S}_m$, then $A\cup B \in \mathcal{S}_m$ and $\overline{A} \in \mathcal{S}_m$.
\item The structure is closed under taking Cartesian products: if $A\in \mathcal{S}_m$ and $B\in \mathcal{S}_n$, then $A\times B\in \mathcal{S}_{m+n}$.
\item Equalities are definable in the structure. Namely, the set $\{(x_1,\ldots,x_m) \mid x_i = x_j\}$ belongs to $\mathcal{S}_m$ for each $i,j\in\{1,\ldots,m\}$.
\item The structure is closed under projections: if $\pi:R^{m+1}\to R^m$ denotes the projection on the first $m$ coordinates, then for all $A\in \mathcal{S}_{m+1}$ the set
$$\pi(A) = \{\underline{x}\in\R^m \mid \exists y\in\R,\ (\underline{x},y)\in A\}$$
belongs to $\mathcal{S}_m$.
\end{enumerate}

The last axiom is called \term{quantifier elimination}.
\end{definition}

For a map $f\colon R^m \to R^n$ we say that it is \term{definable} if its graph is definable:

\[ \{ (\underline{x},\underline{y}) \in R^m\times R^n \mid \underline{y} = f (\underline{x}) \} \in \mathcal{S}_{m+n}. \]

Let us remark that, as was described in the previous sections, Pila's methods involve taking intersections and projections of the initial sets; this shows that definable sets in a structure are a convenient setting for possible generalizations of Pila's results, allowing us to keep track of the regularity of the sets we are working with while we perform operations on them. For our purpose we are particularly interested in the case where the underlying set is $\R$, and we would like to work with structures over $\R$ that are compatible with its properties as an ordered ring. More precisely, we are going to consider structures where 
\begin{itemize}\item The operations $+,\cdot:\R^2\to\R$ are definable. 
\item The singletons $\{x\}$ for all $x\in\R$ are definable.
\item The relation $<$ is definable, in the sense that
for all $m$ and for all $i,j\in\{1,\ldots,m\}$ the set $\{(x_1,\ldots,x_m)\mid x_i< x_j\}$ is definable.
\end{itemize} 
From now on we will understand by ``structure over $\R$'' a structure over the set $R = \R$ satisfying these additional properties. Note that under these assumptions, all semi-algebraic sets, that is, sets given by equations and inequalities involving polynomials, are definable. It is natural to ask whether semi-algebraic sets themselves form a structure over $\R$. This is true, the non-trivial part being axiom $(S4)$, a theorem by Tarski and Seidenberg. 

All semi-algebraic subsets of $\R$ are of one of the following types:
\begin{itemize}\item $\varnothing$;
\item open intervals $(a,b)$ with $a,b\in \R\cup\{\pm\infty\}$;
\item singletons $\{x\}$ for $x\in\R$;
\item finite unions of the previous sets.
\end{itemize}

\begin{definition}A structure over $\R$ is called \term{o-minimal} if the only definable sets in $\R$ are those listed above.
\end{definition}
Though this definition seems to concern only sets in dimension $1$, it imposes strong regularity conditions on definable sets in higher dimensions; for instance, every definable set is a finite union of smooth manifolds,  see \cite{dries}.  On the other hand, the following examples show that the class of o-minimal structures is very rich and allows numerous strong extensions of Pila's results.

\paragraph{Examples:}
					\begin{enumerate}
					\item The simplest (and smallest) example of an o-minimal structure is the structure $\R_{\text{alg}}$ of semi-algebraic sets.

					\item Real semi-analytic sets do not form a structure, because a projection of an analytic set is not necessarily analytic, see \cite{BM88}. Call a set \term{subanalytic} if it is a projection of a semi-analytic set, and \term{globally  subanalytic} if it is a restriction to $\R^m$ of a subanalytic set inside $\mathbb{P}\R^m$. Van den Dries~\cite{Dr86} was, probably, the first one to observe that globally subanalytic sets  form an o-minimal structure, denoted by $\R_\an$. 

					\item The structures $\R_{\text{alg,exp}}$ and $\R_{\text{an,exp}}$ obtained by extending $\R_{\text{alg}}$ and $\R_{\text{an}}$ to make the function $\exp$ definable are o-minimal, as shown by Wilkie, Van den Dries and others in the 1990's, see \cite{DMM94,DM96} for precise references.

					\item
					Peterzil and Starchenko~\cite{PS04} showed definability of the Weierstrass $\wp$-function. Precisely, denote by $\Delta$ the standard fundamental domain for the action of $\text{SL}_2(\Z)$ on the upper half plane $\HH$  (see Subsection~\ref{ssmodu} for more details) and by  $\wp(z,\tau)$ the Weierstrass $\wp$-function of $\Lm = \Z\tau + \Z$.  Then the function $f(x,y,\tau)=\wp(x\tau+y,\tau)$ on the set $[0,1]\times[0,1]\times\Delta$ is definable in $\R_{\text{an},\exp}$. It follows that the $j$-invariant function restricted to $\Delta$ is definable in $\R_{\text{an,exp}}$ as well. 

						\end{enumerate}

\subsection{The theorem of Pila and Wilkie}
 We can now state the principal theorem of Pila and Wilkie from~\cite{PW06}.
\begin{theorem} \label{pilaw}(Pila, Wilkie, 2006) Let a subset $X\subset \R^m$ be definable in some o-minimal structure. Then for every $\eps>0$ there is a constant $c(X,\eps)$ such that
$$|\xtr(\Q,T)|\le c(X,\eps)\,T^{\eps}.$$
\end{theorem}
Applying this to the  o-minimal structure $\R_\an$,  we see that Theorem \ref{higher} we stated earlier is a special case of this one. 

In the sequel we shall often omit for brevity the reference to o-minimality, but will assume it tacitly. By ``definable'' we shall always mean ``definable in some o-minimal structure''.

The main new tool in the proof is an o-minimal version of Gromov's and Yomdin's re-parametrization lemma. 

\begin{definition}Let $X\subset\R^m$ be definable. A partial $\CC^r$-parametrization of $X$ is a definable function $f:(0,1)^m\to X$ which is injective and $\CC^r$. It is said to be \term{bounded} if all its derivatives of order $\le r$ are bounded. A finite set $\{f_1,\ldots,f_k\}$ of (bounded) partial parametrizations of $X$ is called a (\term{bounded}) \term{$\CC^r$-parametrization} of $X$ if the union of the images of all the $f_i$ covers all of $X$. 
\end{definition}

Existence of a $\CC^r$-parametrization (without boundedness condition) is not difficult to establish,  see~\cite{dries}, and the principal hassle is boundedness. Gromov~\cite{Gr87}, improving on a result  of Yomdin~\cite{Yo87}, showed that a semi-algebraic set admits a bounded $\CC^r$-parametrization. Pila and Wilkie \cite{PW06} extended this further to arbitrary definable sets. 

\begin{lemma}(Yomdin--Gromov re-parametrization lemma) A compact definable set has a bounded $\CC^r$-parametrization.
\end{lemma}

Moreover, a similar statement holds for a definable family of definable sets. 

Thanks to the bounded parametrization, one no longer has to struggle with non-uniformity, as we did in Section~\ref{shigher}: all the derivatives are bounded, and one can proceed with the inductive argument almost straightforwardly. 
This is one of the best examples I know when placing the problem in the correct general context  greatly clarifies it, and the proof of a much more general statement appears to be much simpler than that of particular cases.

Later Pila in \cite{Pi09} extended Theorem \ref{pilaw} to algebraic points of bounded degree on definable sets. Denote by $\xtr(\Q,d,T)$ the set of points on $X$ of height\footnote{I intentionally do not specify what I mean here by ``height''; each reader can use her/his favorite definition, the result will always be the same.} bounded by~$T$ and degree (over~$\Q$) bounded by~$d$.

\begin{theorem}\label{pilaalg}(Pila, 2009) Let $d\ge 1$ be an integer and let $X\subset\R^m$ be a definable set. Then for all $\eps>0$ there is a constant $c(X,\eps,d)>0$ such that
$$|\xtr(\Q,d,T)|\leq c(X,\eps,d)\,T^{\eps}.$$
\end{theorem}

The idea is to reduce this to $\Q$-points in a higher dimensional set, precisely, in the symmetric product of $d$ copies of $X$. If~$P$ is a point of degree~$d$ over~$\Q$ and $P_1, \ldots, P_d$ are its conjugates, then the point $(P_1, \ldots, P_d)\in X^d$ gives rise to a $\Q$-rational point on the symmetric product.

\section{The André--Oort conjecture}
Before discussing another very important application of Pila's results, we need to recall some facts on modular curves and complex multiplication, the main reference being \cite{La73}.

\subsection{Modular curves}
\label{ssmodu}

Let $E$ be an elliptic curve over $\C$. Then there is a complex analytic isomorphism $$E(\C)\cong\C/\Lm,$$ where $\Lm$ is some lattice in $\C$. Two elliptic curves $E_1$ and $E_2$ are isomorphic if and only if the corresponding lattices $\Lm_1,\Lm_2$ are equivalent, that is, if there exists $\alpha\in\C^\times$ such that $\Lm_2 = \alpha\Lm_1$. Thus up to isomorphism we can always write $E(\C) =\C/\Lm$ where  $\Lm$ is the lattice $\langle 1,\tau\rangle$ generated by $1$ and an element~$\tau$ of the upper half-plane $\mathcal{H}$. There is a transitive action of $\Gamma(1):=\slz$ on $\mathcal{H}$ given by
$$\left(\begin{array}{cc} a & b \\
												c & d 
												\end{array}\right)\cdot\tau = \frac{a\tau + b}{c\tau + d},$$
and  the modular curve $Y(1) = \Gamma(1)\backslash\HH$ of $\Gamma(1)$-orbits parametrizes isomorphism classes of elliptic curves. One can also describe it more explicitly choosing the following fundamental domain:
$$\Delta = \left\{\tau\in\HH,\ -\frac12\leq\text{Re}(\tau)<\frac12,\ |\tau|>1\right\}\cup\left\{\tau\in\HH,\ |\tau|= 1, -\frac12\le\text{Re}(\tau)\le 0\right\}.$$
Two elliptic curves are isomorphic over $\C$ if and only if they have the same $j$-invariant, and every complex number is the $j$-invariant of some elliptic curve over $\C$. Therefore, the $j$-invariant defines an analytic isomorphism $\Gamma(1)\backslash\HH\cong\C$. 

Define the \term{congruence subgroup}
$$\Gamma_0(N) = \left\{\gamma\in\slz,\ \gamma\equiv \left(\begin{array}{cc} * & * \\
												0 & * 
												\end{array}\right)\pmod N\right\}.$$
												The curve $Y_0(N) = \Gamma_0(N)\backslash\HH$ parametrizes pairs $(E,f)$ where $E$ is an elliptic curve and $f$ is a cyclic $N$-isogeny, that is, an isogeny $f:E\to E'$ (where $E'$ is some other elliptic curve) of degree~$N$ with kernel cyclic of order $N$. We can define a function on $Y_0(N)$ by $j_N(\tau) = j(\tau/N)$. There is an explicit description of $Y_0(N)$ as a plane curve over $\Q$ given by an equation $\Phi_N(x,y) = 0$, where $\Phi_N$ is an irreducible polynomial in $\Z[x,y]$ such that $\Phi_N(j,j_N) = 0$. 

\subsection{Complex multiplication} 
Let $E = E_{\tau}$ be an elliptic curve over $\C$ with corresponding lattice $\Lm = \langle 1,\tau\rangle$. Its endomorphism ring $\mathrm{End}(E)$ contains the ring $\Z$ of rational integers, corresponding to the multiplication-by-$n$ maps for every $n\in\Z$. If $\mathrm{End}(E)$ is strictly bigger than $\Z$, the curve $E$ is said to have \term{complex multiplication}, or to be a \term{CM-curve}.  In this case $\tau$ is necessarily imaginary quadratic, and $\mathrm{End}(E)$ is an order $\mathcal{O}_{\tau}$ in the imaginary quadratic field $K = \Q(\tau)$. In particular, there is a positive integer $f$, called the \term{conductor} of $\tau$, such that $\mathcal{O}_{\tau} = \Z + f\mathcal{O}_K$ where $\mathcal{O}_K$ is the ring of integers of $K$. Writing $\Q(\tau) = \Q(\sqrt{d})$ where $d<0$ is an integer such that $-d$ is square-free, $D = f^2d$ is the discriminant of the order $\mathcal{O}_{\tau}$, called the \term{discriminant} of $\tau$. The following is a fundamental result of the theory of complex multiplication:
\begin{theorem} \label{cm}Assume $E$ has complex multiplication by the order $\mathcal{O}_{\tau}$ of conductor $f$. Then $j(\tau)$ is an algebraic integer of degree $h_{f,K}$, the class number of the order $\mathcal{O}_{\tau}$ of $K$.
\end{theorem}
Moreover, $h_{f,K}$ is related in the following way to the class number $h_K$ of $K$:
\begin{equation}\label{classnb}h_{f,K} = c\cdot h_K \, f\,\prod_{p|f}\left(1 - \left(\frac{d}{p}\right)\frac{1}{p}\right), \qquad c\in \{1,1/2,1/3\}.
\end{equation}
See \cite[Section~7]{Co89} for more details.

\subsection{The André--Oort conjecture}
\begin{definition} A point ${P\in\C^2}$ is said to be a \term{CM-point} if ${P = (j(\tau_1),j(\tau_2))}$ where $\tau_1$ and $\tau_2$ are imaginary quadratic.
\end{definition}
Lines of the form $\C\times\{j(\tau)\}$ or $\{j(\tau)\}\times\C$ contain infinitely many CM-points. So does a modular curve $Y_0(N)$, as we have seen it contains all points of the form $\left(j(\tau),j\left({\tau}/{N}\right)\right).$ André's theorem states that all curves inside $\C^2$ containing infinitely many CM-points must be of one of these types:
\begin{theorem}\label{andre}(André, 1998) If $X\subset\C^2$ is an irreducible curve which is not a horizontal or vertical line nor $Y_0(N)$ for some $N$, then $X$ has finitely many CM-points.
\end{theorem}
The general André--Oort conjecture deals with Shimura varieties, an important special case being products of modular curves $Y_0(N_1)\times\ldots\times Y_0(N_k)$: in this setting it was proved by Edixhoven in 1998 assuming the GRH for imaginary quadratic fields. We are going to concentrate on the Shimura variety $\C^k$ obtained when $N_1 = \ldots = N_k = 1$. 
\begin{definition} \begin{enumerate}\item A \term{special point} in $\C^k$ is a point of the form $(j(\tau_1),\ldots, j(\tau_k))$ where $\tau_i$'s are imaginary quadratic.
\item A \term{special subvariety} of $\C^k$ is a subvariety given by an equation of one of the following forms:
\begin{itemize}\item $x_i = j(\tau)$ for some $i$, where $\tau$ is imaginary quadratic;
\item $\Phi_N(x_i,x_j) = 0$ for some $i,j$.
\end{itemize}
\end{enumerate}
\end{definition}
Special points and special subvarieties are an analogue of torsion points and torsion subvarieties in the statement of the Manin--Mumford conjecture. Now we can state the André--Oort conjecture for $\C^k$, proved by Pila in \cite{Pi11} using his previous results on rational points in definable sets.
\begin{theorem}(André--Oort conjecture) Let $X$ be a subvariety of $\C^k$. Then it has finitely many maximal special subvarieties.
\end{theorem} 
We are going to sketch Pila's proof for $k=2$, when it boils down to proving Theorem \ref{andre}. The ideas are very similar to the ones used for the proof of Manin--Mumford. Consider a map
$$\begin{array}{rcl}\HH^2&\stackrel{\pi}{\longrightarrow}&\C^2.\\
                    (\tau_1,\tau_2)&\mapsto& (j(\tau_1),j(\tau_2))
\end{array}$$
For a non-special irreducible $X\subset\C^2$, put $\tilX = \pi^{-1}(X)\cap(\Delta\times\Delta)$. Any CM-point on $X$ comes from some point of $\tilX$, so it suffices to bound the number of points from $\tilX$ mapping to CM-points in $X$. Note that for any CM-point $P$ on $X$, a point $\widetilde{P}\in\tilX$ mapping to $P$ is of the form $(\tau_1,\tau_2)$ where $\tau_1,\tau_2$ are imaginary quadratic, so $[\Q(\widetilde{P}):\Q]\leq 4.$ We therefore need a bound on the size of $\tilX(\Q,4,T)$ in terms of $T$, and so we must check that we can apply Theorem \ref{pilaalg}. The following is a consequence of the already mentioned theorem by Peterzil and Starchenko about the definability of the Weierstrass $\wp$-function in $\R_{\text{an},\exp}$:
\paragraph{Claim 1:}$\tilX$ is definable.

\bigskip

Moreover, in the same manner as in the proof of Manin--Mumford (an algebraic curve that is not special cannot be $\text{SL}_2(\Z)$-invariant) we obtain

\paragraph{Claim 2:}$\tilX^{\text{tr}} = \tilX$.

\bigskip

Assuming these claims, by Theorem \ref{pilaalg} we get
\begin{equation}\label{upper}|\tilX(\Q,4,T)|\le c(X,\eps)\,T^{\eps}.\end{equation}
Now we are going to bound the number of special points on $X$ from below. For this we can assume that $X$ is defined over some number field $L$. If $P\in X$ is special, all its conjugates over $L$ are special as well. 
Bearing in mind that $P$ is a CM-point, $P = (j(\tau_1),j(\tau_2))$ where we can choose $\tau_1$ and $\tau_2$ inside the fundamental domain $\Delta$, so that $\widetilde{P} = (\tau_1,\tau_2)\in\tilX$. Therefore 
\begin{equation}\label{degrees}[L(P):L]\geq \max \{[\Q(j(\tau_1)):\Q],[\Q(j(\tau_2)):\Q]\}.\end{equation}
According to Theorem \ref{cm}, for $i = 1,2$, 
$$[\Q(j(\tau_i)):\Q] = h_{f_i,\Q(\tau_i)},$$
 where $f_i$ is the conductor of $\tau_i$. The following theorem tells us that the class number of an imaginary quadratic field cannot be too small:
\begin{theorem}(Siegel) Let $K$ be an imaginary quadratic field and fix $\eps'>0$. Then there is a constant ${c = c(\eps')>0}$ such that 
$$h_K\geq c\,|D_K|^{\frac{1}{2}-\eps'},$$
where $D_K$ is the discriminant of the number field $K$.
\end{theorem}
Thanks to formula $(\ref{classnb})$, we then get that 
\begin{equation}\label{degdisc}[\Q(j(\tau_i)):\Q]\geq c_i|D_i|^{\frac{1}{2} - \eps'}\end{equation}
 for some constant ${c_i>0}$, where $D_i = \,f_i^2\,D_{\Q(\tau_i)}$ is the discriminant of $\tau_i$. 
 Moreover, for $\tau$ imaginary quadratic of discriminant~$D$ lying in the fundamental domain, a quick computation shows that the height $H(\tau)$ of $\tau$ is bounded in terms of its discriminant: 
\begin{equation}\label{heightdisc}H(\tau)\leq\sqrt{|D|}.\end{equation}
Putting (\ref{degrees}), (\ref{degdisc}) and (\ref{heightdisc}) together, we get finally for $\eps'$ sufficiently small
$$ [L(P):L]\geq c\max \{H(\tau_1),H(\tau_2)\}^{1-2\eps'}  = c\,H(\widetilde{P})^{1-2\eps'}.$$
If there are infinitely many special points on $X$, the heights $H(\widetilde{P})$ can be arbitrarily large, which with the upper bound $(\ref{upper})$ yields the result.

\bigskip

For further reading on the topic we strongly recommend Scanlon's expository texts \cite{Sc11,Sc12} and Zannier's lecture notes \cite{Za12}. 

As mentioned in the preface, these notes do not claim for any kind of exhaustiveness. For instance, we do not speak at all on Heath-Brown's ``determinant method'', which can be viewed as a ``~$p$-adic'' version of the Bombieri-Pila method; see \cite{HB02,BHS06} and the subsequent work of Browning, Salberger and others.  

While Heath-Brown's method is beautiful and powerful, it mainly applies in the algebraic case, and has yet to show its efficiency in the transcendental case, which is the main topic of these notes. Therefore it is left out.


{\footnotesize






\begin{thebibliography}{}
\bibitem{BM88}\textsc{E. Bierstone, P. D. Milman}, Semianalytic and subanalytic sets, \textit{Pub. Math. I.H.E.S.} \textbf{67} (1988), 5--42.

\bibitem{BP89}
\textsc{E. Bombieri, J. Pila}, The number of integral points on arcs and ovals, \textit{Duke Math. J.} \textbf{59} (1989), 337--357. 

\bibitem{BHS06}
\textsc{T.~D. Browning, D.~R. Heath-Brown,  P. Salberger}, 
Counting Rational Points on Algebraic Varieties,
\textit{Duke Math. J.} \textbf{132} (2006), 545--578.

\bibitem{Co89}
\textsc{D.~A.~Cox}, \textit{Primes of the form ${x^2+ny^2}$},
Wiley, NY, 1989. 

\bibitem{Dr86}\textsc{L. van den Dries},
A generalization of the Tarski-Seidenberg theorem,
and some nondefinability results,
\textit{Bull. Amer. Math. Soc. (N.~S.)} \textbf{15} (1986), 189--193. 

\bibitem{dries}\textsc{L. van den Dries}, \textit{Tame topology and o-minimal structures}, Cambridge Univ. Press, New York, 1998.

\bibitem{DMM94}\textsc{L. van den Dries
A. Macintyre,  D. Marker}, 
The elementary theory of
restricted analytic fields with exponentiation, \textit{Ann. of Math. (2)} \textbf{140} (1994),
183--205.

\bibitem{DM96} \textsc{L. van den Dries, C. Miller}, Geometric categories and o-minimal structures,
\textit{Duke Math. J.} \textbf{84} (1996), 497--540.


\bibitem{Doe27} \textsc{K. Dörge},  Einfacher Beweis des Hilbertschen Irreduzibilitätssatzes,  \textit{Math. Ann.} \textbf{96} (1927),  176--182. 

\bibitem{Ga68} \textsc{A. Gabrièlov}, Projections of semianalytic sets (Russian),  \textit{Funkcional. Anal. i Prilo\v zen.} \textbf{2} (1968),  18--30; \textit{Functional Anal. Appl.} \textbf{2} (1968), 282--291. 

\bibitem{Gr87} \textsc{M. Gromov},  Entropy, homology and semialgebraic geometry, Séminaire Bourbaki 1985/86, no. 663.,  \textit{Astérisque} \textbf{145--146}
(1987), 225--240. 


\bibitem{HW08}
\textsc{G. H.~Hardy,  E. M.~Wright,}
\textit{An introduction to the theory of numbers,}
sixth edition (revised by D.~R.~Heath-Brown and J. H. Silverman, with a foreword by A.~Wiles), Oxford University Press, Oxford, 2008. 

\bibitem{HB02}
\textsc{D.~R.~Heath-Brown}, The density of rational points on curves and surfaces,
\textit{Ann. Math. (2)} \textbf{155} (2002), 553--595.

\bibitem{Ja26} \textsc{V. Jarník,} Über die Gitterpunkte auf konvexen Kurven,  \textit{Math. Z.} \textbf{24} (1926), 500--518. 

\bibitem{La73}\textsc{S. Lang}, \textit{Elliptic functions}, Addison-Wesley, 1973. 

\bibitem{La83} \textsc{S. Lang}, \textit{Fundamentals of Diophantine Geometry}, Springer, 1983. 

\bibitem{La84}
\textsc{M. Laurent}, \'Equations diophantiennes exponentielles, \textit{Invent. Math.} \textbf{78} (1984), 299--327.

\bibitem{Ma84} \textsc{D. W. Masser}, Small values of the quadratic part of the Néron-Tate height on an abelian variety,
\textit{Compositio Math.} \textbf{53} (1984), 153--170. 

\bibitem{PS04}
\textsc{Y.~Peterzil, S. Starchenko,} Uniform definability of the Weierstrass $\wp$ functions and generalized tori of dimension one, \textit{Selecta Math. (N.S.)} \textbf{10} (2004),  525--550. 

\bibitem{Pe06} \textsc{F. Petrov},  On the number of rational points on a strictly convex curve, \textit{Funktsional. Anal. i Prilozhen.} \textbf{40} (2006), 30--42, 95; 
Funct. Anal. Appl. \textbf{40} (2006),  24--33.
\bibitem{Pi05}
\textsc{J. Pila}, Rational points on a subanalytic surface, \textit{Ann. Inst. Fourier} \textbf{55} (2005), 1501--1516.



\bibitem{Pi09}\textsc{J. Pila}, On the algebraic points of a definable set, \textit{Selecta Math. N. S.} \textbf{15} (2009), 151--170.
\bibitem{Pi11} \textsc{J. Pila}, O-minimality and the Andre--Oort conjecture for $\C^n$, \textit{Ann. Math. (2)} \textbf{173} (2011), 1779--1840.
\bibitem{PW06} \textsc{J. Pila, A.J. Wilkie}, The rational points of a definable set, \textit{Duke Math. J.} \textbf{133} (2006), 591--616.

\bibitem{Pl99}
\textsc{A. Plagne,} A uniform version of Jarn\'ik's theorem, \textit{Acta Arith.} \textbf{87} (1999), 255--267.


\bibitem{Ra83}  \textsc{M. Raynaud}, Courbes sur une variété abélienne et points de torsion,  \textit{Invent. Math.} \textbf{71} (1983), 207--233.

\bibitem{Ra83a}  \textsc{M. Raynaud}, Sous-variétés d'une variété abélienne et points de torsion, Arithmetic and geometry, Vol. I, 327--352, Progr. Math. 35, Birkhäuser Boston, Boston, MA, 1983. 

\bibitem{SA94}
\textsc{P.~Sarnak,  S.~Adams}, Betti numbers of congruence groups (with an appendix by Z. Rudnik),
\textit{Israel J. Math.} \textbf{88} (1994), 31--72.

\bibitem{Sc11}
\textsc{Th. Scanlon,}
A Proof of the André--Oort Conjecture Via Mathematical Logic
(after Pila, Wilkie and Zannier), Séminaire Bourbaki,  Avril 2011,
63ème année, 2010-2011, no 1037. 


\bibitem{Sc12}
\textsc{Th. Scanlon,}  Counting special points: Logic, diophantine geometry, and transcendence theory. Bull. Amer. Math. Soc. 49 (2012) 51-71. 



\bibitem{Yo87} \textsc{Y. Yomdin},  $C^k$-resolution of semialgebraic mappings, \textit{Israel J. Math.} \textsl{57} (1987),
301--317.


\bibitem{Za12}
\textsc{U.~Zannier}, \textit{Some Problems of Unlikely Intersections in Arithmetic and Geometry} (with appendices by D. Masser), Ann. Math. Studies 182, Princeton, 2012. 


\bibitem{vieux}\url{http://vieuxgirondin.wordpress.com/2011/07/29/a-mean-value-theorem/}

\end{thebibliography} }
\end{document}